\input amstex
\documentstyle{amsppt}
\magnification=\magstep1
\hoffset=.25truein
\vsize=8.75truein
\NoBlackBoxes
\nologo

\def\bs{{\backslash}}

\def\IN{{\Bbb N}}
\def\IR{{\Bbb R}}

\def\cC{{\Cal C}}
\def\cP{{\Cal P}}
\def\cT{{\Cal T}}
\def\cF{{\Cal F}}

\def\undm{{\underline m}}
\def\undtau{{\underline{\tau}}}
\def\undcT{{\underline{\cT}}}
\def\barE{{\overline E}}
\def\hatE{{\widehat E}}
\def\tilE{{\widetilde E}}
\def\endprf{{\hfill $\square$}}

\def\rint{{\text{int}\, }}
\def\Span{{\text{Span}\, }}

\topmatter
\title{Differentiable functions defined in closed sets.  A problem
of Whitney}\endtitle
\author{Edward Bierstone, Pierre D. Milman and 
Wies\l aw Paw\l ucki}\endauthor
\rightheadtext{Bierstone, Milman, Paw\l ucki}
\leftheadtext{Differentiable functions defined in closed sets}
\endtopmatter

\address{E.B. and P.M.: Department of Mathematics, University of
Toronto, Toronto, Ontario, Canada M5S 3G3}\endaddress
\email{bierston\@math.toronto.edu, milman\@math.toronto.edu}\endemail

\address{W.P.: Uniwersitet Jagiello\'nski, Instytut Matematyki,
ul. Reymonta 4, 30059 Krak\'ow, Poland}\endaddress
\email{Wieslaw.Pawlucki\@im.uj.edu.pl}\endemail

\thanks
Research partially supported by the following
grants: E.B. -- NSERC OGP0009070, P.M. -- NSERC OGP0008949 and 
the Killam Foundation, W.P. -- KBN 5 PO3A 005 21.
\endthanks


\abstract
In 1934, Whitney raised the question of how to
recognize whether a function $f$ defined on a closed subset
$X$ of $\IR^n$ is the restriction of a function of class 
$\cC^p$. A necessary and sufficient criterion was given in
the case $n=1$ by Whitney, using limits of finite differences,
and in the case $p=1$ by Glaeser (1958), using limits of secants. 
We introduce a necessary geometric criterion, for general
$n$ and $p$, involving limits
of finite differences, that we conjecture is sufficient at
least if $X$ has a ``tame topology''. We prove that, if 
$X$ is a compact subanalytic set, then there exists 
$q = q_X(p)$ such that the criterion of order $q$ implies
that $f$ is $\cC^p$. The result gives a new approach to
higher-order tangent bundles (or bundles of differentiable
operators) on singular spaces.
\endabstract

\endtopmatter
\document
\baselineskip=18pt 

\head 1. Introduction\endhead

In 1934, Hassler Whitney published three pioneering
articles on criteria for a function $f:X\to \IR$, where 
$X$ is a closed subset of $\IR^n$, to be the restriction 
of a function of class $\cC^p$, \cite{W1}, \cite{W2}, \cite{W3}. 
($\cC^p$ means continuously differentiable to order $p$, where 
$p \in \IN$.) Whitney's extension theorem \cite{W1} 
gives a necessary and sufficient condition for a field of 
polynomials $\sum_{|\alpha| \le p} f_\alpha (a) 
(x-a)^\alpha $, $a\in X$, where $f_0 = f$, to be the field of 
Taylor polynomials of a $\cC^p$ function.  (We use multiindex 
notation: $\alpha = (\alpha_1,\dots,\alpha_n) \in \IN^n$, 
$|\alpha | = \alpha_1 + \dots + \alpha_n $ and 
$x^\alpha = x^{\alpha_1}_1 \cdots x^{\alpha_n}_n $.
$\IN$ denotes the nonnegative integers.) In general,
the functions $f_\alpha$
are, of course, not uniquely determined by $f$.
In {\it Differentiable functions defined in closed sets. I}
\cite{W2}, Whitney raises the deeper question of a necessary 
and sufficient criterion involving only the values of $f$, 
and he answers the question in the case $n=1$. 
Whitney proves that, if $X$ is a closed subset of the real 
line, then $f$ extends to a $\cC^p$ function if and only if 
the limiting values of all $p$'th divided differences 
$[x_0 , x_1,\dots,x_p] f$, where the $x_i \in X$ and 
$x_i \ne x_j$ if $i\ne j$, define a continuous function on the 
diagonal $\{ x_0 = x_1 = \cdots = x_p \}$. 
($[x_0 , x_1,\dots, x_n]f = p! c_p$, where 
$P(x) = c_0 + c_1 x + \cdots + c_p x^p$ is the 
{\it Lagrange interpolating polynomial} for $f$ at the 
points $x_0$, $x_1,\dots , x_p$; i.e., the unique polynomial 
of degree at most $p$ such that $P(x_i) = f(x_i)$, 
$i=0,\dots,p$.)

``Differentiable functions defined in closed sets. II'' never 
appeared, and up to now the only significant progress on Whitney's 
problem following \cite{W3} seems to have been the beautiful 
theorem of Georges Glaeser (\cite{G}, 1958) which solves the 
problem in the case $p=1$ 
(cf. \cite{Br}. See also Remark 2.3.) 
Glaeser defines a ``(linearized) paratangent
bundle'' $\tau (X)$ using limits of secant lines.
(See Section 3 below.) Suppose that $f$ is continuous and 
let $\tau (f)$ denote the paratangent bundle of the graph of 
$f$.  Then $\tau (f)$ can be regarded as a bundle over $X$, 
and $\tau (f) \subset \tau (X) \times \IR$ 
(but $\tau (f)$ does not necessarily project onto 
$\tau (X)$). Glaeser proves that $f$ is the restriction of a 
$\cC^1$ function if and only if $\tau (f)$ defines a function 
$\tau (f) : \tau (X) \to \IR$ (i.e., $\tau (f)$ is the graph 
of a function $\tau(X) \to \IR$; it will be convenient to 
identify a function with its graph). 

In this article, we introduce a ``(linearized) paratangent bundle
of order $p$'' $\tau^p (X)$, for any $p\in \IN$, using 
limits of finitely supported distributions with values in the 
dual space $\cP^*_p$ of the space $\cP_p = \cP_p (\IR^n)$ of 
polynomial functions on $\IR^n$ of degree at most $p$
(Section 4 below).  Each fibre $\tau^p_a (X)$, $a\in X$, 
is a linear subspace of $\cP^*_p$.  Our construction involves 
a new interpretation of the remainder condition in Whitney's
extension theorem. 

To every function $f:X\to \IR$, we associate a bundle 
$\nabla^p f \subset \tau^p (X) \times \IR$. 

\proclaim{Conjecture}
$f$ is the restriction of a $\cC^p$ function if and only if 
$\nabla^p f$ defines a function $\nabla^p f: \tau^p (X) \to 
\IR$.

Moreover, if $\nabla^p f: \tau^p (X) \to
\IR$ and $\nabla^p_a f = 0$, for some $a \in X$, then there
exists $F\in \cC^p (\IR^n)$ such that $F|X = f$ and $T^p_a F = 
0$, where $T^p_a F$ denotes the Taylor polynomial of order $p$
of $F$ at $a$.
\endproclaim

Necessity of the criterion $\nabla^p f: \tau^p (X) \to
\IR$ is not difficult; the following 
theorem is proved in Section 4.

\proclaim{Theorem 1.1}
If $f:X \to \IR$ extends to a $\cC^p$ function, then 
$$ \nabla^p f : \ \tau^p (X) \to \IR\ . $$
Moreover, if $F\in \cC^p (\IR^n)$ and $F|X = f$, then, 
for all $a\in X$ and $\xi \in \tau^p_a (X) \subset \cP^*_p$, 
$$ \nabla^p f (\xi) \ = \ \xi (T^p_a F)\ . $$
\endproclaim 

The converse direction is true if $X$ is a $\cC^p$ submanifold.
In Section 4, we prove more precisely: 

\proclaim{Theorem 1.2}
If $X\subset M$, where $M$ is a $\cC^p$ submanifold of $\IR^n$ and 
$X$ is the closure of its interior in $M$, 
then $f: X\to \IR$ extends to a 
$\cC^p$ function if and only if $\nabla^p f : \tau^p (X) \to \IR$.
\endproclaim

In order to make the conjecture tractable in general, it is reasonable
to restrict to closed sets $X$ that have a ``tame'' geometry
(``g\'eom\'etrie mod\'er\'ee''); for example, closed subanalytic 
sets or, more generally, closed sets that are definable in an 
$o$-minimal structure (cf. \cite{vdD}).  
Our main result is the following theorem (proved in Section 5). 

\proclaim{Theorem 1.3}
Let $X$ be a compact subanalytic subset of $\IR^n$.  Then there
is a function $q = q_X(p) \geq p$ from $\IN$ to itself such that, 
if $f: X \to \IR$, \ $q \geq q_X(p)$ and 
$$ \nabla^q f : \ \tau^q (X) \to \IR\ , $$ 
then $f$ extends to a $\cC^p$ function. 
If, moreover, $a \in X$ and $\nabla^q_a f = 0$, then 
$f$ extends to a $\cC^p$ function that is $p$-flat at $a$.
\endproclaim 

The novelty of Theorem 1.3 lies in the construction of 
$\tau^p (X)$ and $\nabla^p f$.  Let $X$ be a compact 
subanalytic subset of $\IR^n$.  Then there is a compact 
real analytic manifold $M$ such that $\dim M = \dim X$, and a 
real analytic mapping $\varphi : M \to \IR^n$ such that 
$\varphi (M) = X$ (by the Uniformization Theorem \cite{BM1, Thm. 0.1}).
Let $g = f\circ \varphi$.  We prove that if $\nabla^p f : 
\tau^p (X) \to \IR$, then: 

(1)\  $\nabla^p g : \tau^p (M) \to \IR$ (Theorem 5.2);
therefore, $g\in \cC^p (M) $ by Theorem 1.2.  (All notions 
make sense for manifolds.)

(2)\ $g$ is {\it formally a composite
with} $\varphi$; i.e., for all $a\in X$, there exists 
$P\in \cP_p (\IR^n)$ such that $g-P \circ \varphi$ is $p$-flat 
at every point $b\in \varphi^{-1} (a)$ (Corollary 5.3).

Theorem 1.3 is then a consequence of the following composite 
function theorem \cite{BMP}: 
There is a function $q = q_\varphi (p)$ such that if $g\in \cC^q (M)$
is formally a composite with $\varphi$, then there exists
$F \in\cC^p (\IR^n)$ such that $g = F \circ \varphi$. (Moreover,
if $S$ is a finite subset of $X$ and
$g$ is $q$-flat on $\varphi^{-1} (S)$, then there exists $F$ with
the additional property that $F$ is $p$-flat on $S$.)

Let $\cC^{(\infty)} (X) = \bigcap_{p\in \IN} \cC^p (X)$, 
where $\cC^p (X)$ denotes the space of restrictions of $\cC^p$
functions to $X$. 

\proclaim{Corollary 1.4}
If $X\subset \IR^n$ is a closed subanalytic set and $f: X\to \IR$, 
then $f\in \cC^{(\infty)} (X)$ if and only if 
$\nabla^p f : \tau^p (X) \to \IR$, for all $p\in \IN$. 
\endproclaim 

Of course $\cC^\infty (X) \subset \cC^{(\infty)} (X)$, 
where $\cC^\infty (X)$ denotes the restrictions of $\cC^\infty$
functions to $X$, and $\cC^{(\infty)} (X) = \cC^\infty (X)$
if $X$ is a $\cC^\infty$ submanifold (but not in general \cite{P2}).  
Among closed subanalytic 
sets, equality characterizes the proper subclass of sets that have
a ``semicoherent'' (or stratified coherent) structure 
\cite{BM2}; see Remarks 2.7 below. In Section 2, 
we use Theorem 1.3 to compare the
paratangent bundle $\tau^p (X)$ with another natural idea of
a bundle of differential operators on a singular space. 
Surprisingly, uniformity of a ``Chevalley estimate'' (which also
characterizes the class of semicoherent subanalytic sets) is
related to important stability properties of these bundles.
(See Remarks 2.7, Theorem 2.9 and Corollary 2.10.)

The space $\cC^{(\infty)} (X)$ seems to be an interesting
function space for a closed set $X$ that is definable in an 
$o$-minimal structure.  A definable set has a $\cC^p$ cell 
decomposition for every $p$, but $\cC^\infty$ cell decomposition, 
in general, is unknown and likely untrue. 
Theorem 1.3 provides strong evidence for the conjecture
above in the case of a definable set $X$.  (See Final Remarks 5.5.)
The loss of differentiability in Theorem 1.3 is related to the 
use of \cite{BMP} via a uniformization of $X$.

\head 2. Geometric and algebraic paratangent bundles \endhead

In this section, we introduce a ``Zariski paratangent bundle''
$\cT^p (X)$ -- a (higher-order) $\cC^p$ analogue of the Zariski 
tangent bundle studied in algebraic geometry -- and we use the
results above to compare the paratangent bundle $\tau^p (X)$
with $\cT^p (X)$. (Whitney \cite{W4} makes such a comparison in
order $1$, for various notions of tangent spaces to an analytic
variety.) It seems less interesting to use the Zariski paratangent
bundle to provide a criterion to recognize whether a function
$f:X \to \IR$ is the restriction of a $\cC^p$ function because
$\cT^p (X)$ is defined already in terms of the ideal of
$\cC^p$ functions vanishing on $X$ (hence essentially in terms of
the space of restrictions to $X$ of $\cC^p$ functions); see
Remark 2.3. The interest is rather in the opposite direction -- to
use the conjecture or the results in Section 1, involving limits
of finite differences, to get a better understanding of
higher-order tangent bundles (or bundles of differential 
operators) on singular spaces. This section is not used in the
rest of the paper, except in Remarks 4.13(1). 

We use the notation of Section 1.
If $V \subset \cP_p = \cP_p (\IR^n)$, let $V^\perp$ denote the
orthogonal complement of $V$ in the dual space $\cP_p^*$. Let
$X$ be a closed subset of $\IR^n$. Let $I^p(X) \subset \cC^p 
(\IR^n)$ denote the ideal of $\cC^p$ functions that vanish on
$X$. 

\definition{Definition 2.1}
The {\it Zariski paratangent bundle} of {\it order}
$p$, $\cT^p (X)$, is the subbundle of $X \times \cP_p^*$  with
fibre $\cT^p_a (X) = (T^p_a I^p(X))^\perp$, for each $a \in X$.
(See Definition 3.1 below.) 
\enddefinition 

The bundle $\cT^p (X)$ is closed in 
$X \times \cP_p^*$ because, for all $h \in I^p(X)$, $Z(h) :=
\{ (a, \xi) \in X \times \cP_p^* : \ \xi (T^p_a h) = 0 \}$ is
closed, and $\cT^p (X) = \bigcap_{h \in I^p(X)} Z(h)$. Moreover,
if $a \in X$, then
$$
\tau^p_a (X) \ \subset \ (T^p_a I^p(X))^\perp \ \subset \ \cP_p^*
$$
(by Definition 4.12 or Theorem 1.1); 
thus $\tau^p (X) \subset \cT^p (X)$. Both
$\tau^p (\cdot)$ and $\cT^p (\cdot)$ are functors on the
category of closed (or locally closed) subsets of Euclidean
spaces, with morphisms given by the restrictions of $\cC^p$
mappings; cf. Theorem 5.2(1) below.

Consider $q \geq p$. Let $a \in \IR^n$. 
Let $\undm^{p+1}_a \subset \cP_q$ denote the subspace of polynomials
of order at least ${p+1}$ at $a$; i.e.,
$$
\undm^{p+1}_a \ = \ \{ P \in \cP_q : \ (D^\alpha P)(a) = 0, \ 
|\alpha| \leq p \}.
$$
There is a projection $\cP_q \to \cP_p$ defined by truncating
terms in $x-a$ of order $> p$ in the expression of any $P \in 
\cP_q$ as a polynomial in $x-a$. Let $\cP_p^* \hookrightarrow 
\cP_q^*$ denote the embedding dual to this projection.
The projection induces an isomorphism $\cP_q / \undm^{p+1}_a 
\cong \cP_p$, and the embedding $\cP_p^* \hookrightarrow \cP_q^*$
has image $(\undm^{p+1}_a)^\perp$.

Let $\tau^q (X)_p$ denote the subbundle of $X \times \cP_p^*$
with fibre 
$$
\tau^q_a (X)_p \ := \ \tau^q_a (X) \cap (\undm^{p+1}_a)^\perp \ ,
$$
for each $a \in X$, where $(\undm^{p+1}_a)^\perp$ is identified with 
$\cP_p^*$ via the embedding above. Of course, $\tau^p (X)_p 
= \tau^p (X)$ . 

\proclaim{Lemma 2.2}
$\tau^q (X)_p$ is a
closed subbundle of $X \times \cP_p^*$.
\endproclaim

\demo{Proof}
There is a continuous bundle mapping (cf. Definition 4.23) 
$X \times \cP_p^* \to X \times \cP_q^*$ with closed image,
where, for each $a \in X$, the fibre $\cP_p^*$ over $a$ is
embedded in $\cP_q^*$ as above. Via this
mapping, $X \times \cP_p^*$ is a closed subbundle of 
$X \times \cP_q^*$. Of course, 
$\tau^q (X)_p = \tau^q (X) \cap (X \times \cP_p^*)$.
\endprf
\enddemo

If $F \in \cC^q(\IR^n)$, let $T^p_a F$ denote the Taylor polynomial
of order $p$ of $F$ at $a$. Let $\cT^q (X)_p \subset 
X \times \cP_p^*$ denote the bundle with fibre 
$$
\cT^q_a (X)_p \ := \ (T^p_a I^q(X))^\perp
\ = \ \cT^q_a (X) \cap (\undm^{p+1}_a)^\perp \ ,
$$ 
for each $a \in X$. 
Then $\cT^q (X)_p$ is closed
in $X \times \cP_p^*$, and
$$
\tau^q (X)_p \ \subset \ \cT^q (X)_p \ .
$$

\remark{Remark 2.3}
It is possible to formulate various criteria for the existence
of a $\cC^p$ extension of $f:X \to \IR$ involving only the
values of $f$ on $X$. The following is essentially tautological.

{\it Implicit function criterion.} Suppose that $f$ is continuous.
Then $f$ is the restriction of a $\cC^p$ function if and only if,
for every $a \in X$, there is a neighbourhood of $(a, f(a))$ in
$\IR^n \times \IR$ in which the graph of $f$ lies in a $\cC^p$
submanifold whose tangent space at $(a, f(a))$ contains no vertical
vector (i.e., no derivation in the vertical direction). 

The result of O'Farrell  and Watson \cite{O'FW} 
is a closely related criterion that
can be expressed in terms of the Zariski paratangent bundle 
$\cT^p (X)$. Suppose that $f$ is continuous. Let $\cT^p (f)$
denote the Zariski paratangent bundle of order $p$ of the graph
of $f$. Then $\cT^p (f)$ can be regarded as a bundle over $X$, and
the projection $\IR^n \times \IR \to \IR^n$ induces a bundle
mapping $\pi: \ \cT^p (f) \to \cT^p (X)$. The theorem of \cite{O'FW}
asserts essentially that $f$ is the restriction of a $\cC^p$ 
function if and only if $\pi$ is bijective. ``Only if'' is obvious.
On the other hand, if $\pi$ is bijective, then, for each $a \in X$,
$\cT^p_a (f)_1$ contains no derivation in the vertical direction.
It follows from Definition 2.1 
that, in some neighbourhood of $(a, f(a))$, the graph of $f$ lies
in a submanifold whose tangent space at $(a, f(a))$ contains no
vertical derivation. The result follows from the
implicit function criterion.
\endremark
\medskip

Our conjecture in Section 1 implies that
$$
\tau^p (X) \ = \ \cT^p (X) \ :
$$
We need only show that $\cT^p (X) \subset \tau^p (X)$. Suppose
that $P \in \tau^p_a (X)^\perp \subset \cP_p$; we have to show 
that $P \in T^p_a I^p(X)$. Let $f = P|X$. The conjecture asserts
that there exists $F\in \cC^p (\IR^n)$ such that $F|X = f$ and 
$T^p_a F = 0$. Let $G = P - F$. Then $G \in I^p(X)$, and $T^p_a G
= P - T^p_a F = P$.

The following is a corollary of Theorem 1.3.

\proclaim{Corollary 2.4}
Suppose that $X$ is a compact subanalytic subset of $\IR^n$.
Let $q = q_X(p)$ denote a function satisfying the assertion
of Theorem 1.3. If $s \geq p$ and $q \geq q_X(s)$, then 
$$
\cT^s (X)_p \ \subset \ \tau^q (X)_p \ .
$$
\endproclaim

\demo{Proof}
Let $a \in X$ and let $P \in \tau^q_a (X)_p^\perp \subset
\cP_q$. We have to show that $P \in T^p_a I^s(X) + \undm^{p+1}_a$. 
Recall that
$\tau^q_a (X)_p = \tau^q_a (X) \cap (\undm^{p+1}_a)^\perp$;
thus $P \in \tau^q_a (X)^\perp + \undm^{p+1}_a$; i.e., there
exists $Q \in \tau^q_a (X)^\perp \subset \cP_q$ such that 
$T^p_a Q - P \in \undm^{p+1}_a$. Let $f = Q|X$. 
By Theorem 1.1, $\nabla^q f : \tau^q (X) \to \IR$
and $\nabla^q_a f = 0$ (since $Q \in \tau^q_a (X)^\perp$). By
Theorem 1.3, there exists $F\in \cC^s (\IR^n)$ such that $F|X = f$ 
and $T^s_a F = 0$. Set $G = Q - F$. Then $G \in I^s(X)$, and 
$T^p_a G - P = T^p_a Q - P \in \undm^{p+1}_a$.
\endprf
\enddemo

Following the viewpoint of Corollary 1.4 above, we can also
introduce
$$
\align
\undtau^p (X) \ &:= \ \bigcup_{q \geq p} \tau^q (X)_p \ ,\\
\undcT^p (X) \ &:= \ \bigcup_{q \geq p} \cT^q (X)_p \ .\\
\endalign
$$
Then $\undtau^p (X) \subset \undcT^p (X)$.

\proclaim{Corollary 2.5}
If $X$ is a closed subanalytic subset of $\IR^n$, then
$\undtau^p (X) = \undcT^p (X)$.
\endproclaim

\definition{Definition 2.6}
Let $\IR [[x-a]]$ denote the ring of formal power series in
$(x_1 - a_1, \ldots, x_n - a_n)$. If $F \in \IR [[x-a]]$, let
$T^p_aF(x)$ denote the Taylor polynomial of order $p$ of $F$
at $a$; i.e., the polynomial of degree $\leq p$ obtained by
truncating the terms of $F$ of order $> p$ in $x-a$. Suppose 
$a \in X$. We define the {\it formal local ideal} $\cF_a(X)$
of $X$ at $a$ as $\{ F \in \IR [[x-a]]: \ T^p_aF(x) =
o(|x-a|^p),$ where $x \in X,$ for all $p \in \IN\}$. (See
\cite{BM2, Lemma 6.1}.)
\enddefinition

\remark{Remarks 2.7}
For each $q \geq p$, $\tau^q (X)_p \subset \tau^{q+1} (X)_p$
(as in Remarks 4.13(1) below) and 
$\cT^q (X)_p \subset \cT^{q+1} (X)_p$; i.e.,
$\{\tau^q (X)_p\}_{q \geq p}$ and
$\{\cT^q (X)_p\}_{q \geq p}$ are increasing sequences of closed
subbundles of $X \times \cP_p^*$. If $\{\tau^q (X)_p\}_{q \geq p}$
(respectively, $\{\cT^q (X)_p\}_{q \geq p}$) stabilizes,
then $\undtau^p (X)$ (respectively,
$\undcT^p (X)$) is closed. Since $\tau^q (X)_p \subset \cT^q (X)_p$,
for all $q \geq p$, it follows from Corollary 2.4 that, if $X$ is
a compact subanalytic set, then, for all $p \in \IN$, 
$\{\tau^q (X)_p\}_{q \geq p}$ stabilizes if and only if
$\{\cT^q (X)_p\}_{q \geq p}$ stabilizes. 

We say that a closed
subanalytic subset $X$ of $\IR^n$ is {\it semicoherent} if it has
a locally finite subanalytic stratification such that the formal
local ideals $\cF_a(X)$ are generated over each stratum by finitely
many subanalytically parametrized formal power series \cite{BM2, 
Definition 1.2}. In Corollary 2.10 below, we show that a compact
subanalytic subset $X$ of $\IR^n$ is semicoherent if and only if
the sequence $\{\cT^q (X)_p\}_{q \geq p}$ (or the sequence
$\{\tau^q (X)_p\}_{q \geq p}$) stabilizes, for all $p \in \IN$.

For compact subanalytic sets $X$ in general, however, 
$\undcT^p (X)$ is not necessarily closed, and
$\{\cT^q (X)_p\}_{q \geq p}$ does not necessarily 
stabilize even if $\undcT^p (X)$ is closed. If $n \leq 4$,
$\dim X \leq 2$, or $\dim X \geq n-1$, 
then $X$ is semicoherent. In $\IR^5$, consider
any sequence of distinct points $\{a_j\}$ tending to the origin along
some line. By the construction of \cite{P1}, there is a compact
$3$-dimensional subanalytic subset $X$ of $\IR^5$ such that $X$
is not semicoherent, $X$ is semicoherent outside $0$, and
$\cF_a(X) = 0$ if and only if $a \in \{a_j\}$ (cf. \cite{BM2,
Examples 1.29}). By Corollary 2.10, 
$\undcT^p (X) \backslash \undcT^p_0 (X)$ is closed in 
$(X \backslash 0) \times \cP_p^*$, for all $p$. For all $p$ and
$j$, $\undcT^p_{a_j} (X) = \cP_p^*$, by Lemma 2.8, since 
$\cF_{a_j}(X) = 0$. It follows that, if $0 \in \{a_j\}$, then
$\undcT^p (X)$ is closed, for all $p$. On the other hand, if
$0 \notin \{a_j\}$, then $\cF_0(X) \neq 0$, so there exists
$p$ such that $\undcT^p_0 (X) \neq \cP_p^*$, and it follows
that $\undcT^p (X)$ is not closed.
\endremark
\medskip

Let $X$ be a closed subset of $\IR^n$ and let $a \in X$. 
Then $\{T^p_a I^q (X)\}_{q \geq p}$ is a decreasing sequence
of linear subspaces of $\cP_p$. Let $s_X(a,p)$ denote the
smallest integer $s \geq p$ such that $T^p_a I^q (X) = 
T^p_a I^s (X)$, for all $q \geq s$. Of course, $s_X(a,p)$ is
the smallest integer $s \geq p$ such that $\cT^s_a (X)_p = 
\bigcup_{q \geq p} \cT^q_a (X)_p$.

\proclaim{Lemma 2.8}
Suppose that $X$ is subanalytic. Let $a \in X$ and let
$s \geq p$. Then $s \geq s_X(a,p)$ if and only if
$$
T^p_a I^s (X) \ = \ T^p_a \cF_{a}(X) \ .
$$
\endproclaim

\demo{Proof}
First we show that, if $s \geq s_X(a,p)$, then 
$T^p_a I^s (X) \subset T^p_a \cF_{a}(X)$, {\it for any}
closed $X \subset \IR^n$. Let $p_0 = p$ and $p_j = 
s_X(a,p_{j-1})$, for all $j \geq 1$. Let $P = T^p_a f_0$,
where $f_0 \in I^{p_1}(X)$. We have to show that $P = T^p_a F$,
where $F \in \cF_{a}(X)$. By the definition of $\{p_j\}$, for
all $j \geq 1$,\ $T^{p_{j-1}}_a I^{p_j}(X) = 
T^{p_{j-1}}_a I^{p_{j+1}}(X)$, so that we can inductively find
$f_j \in I^{p_{j+1}}(X)$, $j \geq 1$, such that $T^{p_{j-1}}_a 
f_{j-1} = T^{p_{j-1}}_a f_j$. Let $Q_j = T^{p_j}_a f_j$,\ $j
\geq 0$, so that $Q_j(x) = o(|x-a|^{p_j})$, where $x \in X$,
and $Q_{j-1} = T^{p_{j-1}}_a Q_j$, \ $j \geq 1$. Take $F \in
\IR [[x-a]]$ such that $Q_j = T^{p_j}_a F$, for all $j$. Then
$P = T^p_a F$ and, for all $j$, $T^{p_j}_a F = Q_j = 
o(|x-a|^{p_j})$, where $x \in X$, as required.

We can assume that $X$ is a compact
subanalytic set. Then $T^q_a \cF_{a}(X) \subset T^q_a I^q(X)$,
for all $q$, as follows: There is a compact real analytic
manifold $M$ and a real analytic mapping $\varphi:
\ M \to \IR^n$ such that $\varphi (M) = X$. Let $F \in
\cF_{a}(X)$. Take $f \in \cC^\infty (\IR^n)$ such that $F$
is the formal Taylor series of $f$ at $a$. Set $g = f \circ 
\varphi$. Then $g$ is flat on ${\varphi}^{-1}(a)$. By
\cite{BMP}, for all $q \in \IN$, there exists $f_q \in
\cC^q(\IR^n)$, such that $T^q_a f_q = 0$ and
$g = f_q \circ \varphi$. Then $f - f_q \in I^q(X)$ and
$T^q_a (f - f_q) = T^q_a f = T^q_a F$.

Thus $T^p_a \cF_{a}(X) \subset T^p_a I^q(X)$,
for all $q \geq p$. If $s \geq s_X(a,p)$, then $T^p_a I^s (X)
\subset T^p_a \cF_{a}(X)$; hence $T^p_a I^s (X) = T^p_a \cF_{a}(X)$.
Conversely, if $T^p_a I^s (X) = T^p_a \cF_{a}(X)$, then, for all
$q \geq s$, $T^p_a I^q (X) \subset T^p_a I^s (X) = T^p_a \cF_{a}(X) 
\subset T^p_a I^q (X)$, so that $T^p_a I^q (X) = T^p_a I^s (X)$,
and hence $s \geq s_X(a,p)$.
\endprf
\enddemo

Let $X$ denote a compact subanalytic subset of $\IR^n$.
Let $\varphi: \ M \to \IR^n$ be a real analytic mapping
from a compact real analytic manifold $M$, 
such that $\varphi (M) = X$.
By \cite{BMP}, for all $p \in \IN$, there exists $q \geq p$
with the following property: if $a \in X$ and $g \in \cC^q (M)$
such that $g$ is formally a composite with $\varphi$ and $g$
is $q$-flat on ${\varphi}^{-1}(a)$, then there exists $f \in
\cC^p(\IR^n)$ such that $g = f \circ \varphi$ and $f$ is $p$-flat
at $a$. Let $q_{\varphi}(p)$ denote the least such $q$.

By a lemma of Chevalley (cf. \cite{BM2, Section 6}), for all $k \in
\IN$, there exists $l \in \IN$, $l \geq k$,
with the following property:
for every polynomial
$F(x)$ such that $F(x) = o(|x-a|^l)$, where
$x \in X$, there exists $G \in \cF_{a}(X)$ such that 
$T^k_a F = T^k_a G$.
Given $k$, let $l_X(a, k)$ denote the least such $l$. We call
$l_X(a, k)$ a {\it Chevalley estimate}.

\proclaim{Theorem 2.9}
Let $X$ be a compact subanalytic subset of $\IR^n$ and let 
$\varphi: \ M \to \IR^n$ be a real analytic mapping 
as above. Then, for
all $a \in X$ and all $p \in \IN$,
$$
s_X(a,p) \ \leq \ l_X(a,p) \ \leq \ q_{\varphi}(s_X(a,p)) \ .
$$
\endproclaim

\demo{Proof}
For the first inequality, let $s = l_X(a,p)$; by Lemma 2.8,
it is enough
to show that if $f \in I^s(X)$, then $T^p_a f \in T^p_a \cF_{a}(X)$.
Let $P = T^s_a f$. Then $P(x) = o(|x-a|^s)$, where $x \in X$, so
the result follows from the definition of $l_X(a,p)$.

For the second inequality, let $s = s_X(a,p)$ and 
$q = q_{\varphi}(s)$. Let $F$ denote a polynomial
such that $F(x) = o(|x-a|^q)$, where $x \in X$. By
Lemma 2.8, it is enough to show that $T^p_a F \in T^p_a I^s(X)$.
Let $g = F \circ \varphi$. Then $g \in \cC^\infty (M)$ and
$g$ is $q$-flat on ${\varphi}^{-1}(a)$. By \cite{BMP}, there
exists $f \in \cC^s(\IR^n)$ such that $T^s_a f = 0$ and 
$g = f \circ \varphi$. Then $f = F$ on $X$ and $T^s_a f = 0$,
so that $T^p_a F = T^p_a (F-f) \in T^p_a I^s(X)$.
\endprf
\enddemo

\proclaim{Corollary 2.10}
Let $X$ be a compact subanalytic subset of $\IR^n$.

(1)\ Let $p \in \IN$. Then the increasing union $\bigcup_{q \geq p}
\cT^q (X)_p$ stabilizes if and only if there exists $l \in \IN$
such that $l_X(a,p) \leq l$, for all $a \in X$.

(2)\ The increasing union $\bigcup_{q \geq p} \cT^q (X)_p$ (or the
increasing union $\bigcup_{q \geq p} \tau^q (X)_p$) stabilizes,
for all $p \in \IN$, if and only if $X$ is semicoherent.
\endproclaim

\demo{Proof}
The first assertion is immediate from Theorem 2.9. By \cite{BM2,
Theorem 1.13}, $X$ is semicoherent if and only if there exists
a {\it uniform Chevalley estimate}; i.e., a function $l = l_X(k)$ such
that $l_X(a,k) \leq l_X(k)$, for all $k \in \IN$ and $a \in X$. 
So the second statement also follows.
\endprf
\enddemo

\head 3. Glaeser's construction\endhead 

Let $X$ be a metric space and let $V$ be a real 
vector space of finite dimension $r$. 

\definition{Definition 3.1}
A {\it bundle} ({\it of linear subspaces of} $V$)
{\it over} $X$ is a subset $E$ of $X\times V$ such that, 
for all $a\in X$, the {\it fibre} $E_a := \{ v\in V : (a,v) \in 
E\}$ is a linear subspace of $V$. 
\enddefinition

\definition{Definition 3.2}
A {\it Glaeser operation} (on bundles of linear subspaces 
of $V$ over $X$) is an operation $\rho$ that associates
to each bundle $E$ a bundle $\rho(E)$ such that: 

(1)\  $\barE \subseteq \rho (E)$; 

(2)\ $\rho$ is {\it local}; i.e., 
if $E$, $F$ are bundles over $X$ and $E_a = F_a$ for all 
$a\in U$, where $U \subset X$ is open, then 
$\rho (E)_a = \rho (F)_a $ for all $a\in U$. 
\enddefinition

We include a proof of the following lemma of Glaeser \cite{G}
because we use it in Sections 4 and 5. 

\proclaim{Lemma 3.3}
Let $E$ be a bundle of linear subspaces of $V$ over $X$, 
and let $\rho$ be a Glaeser operation. 
Write $\rho^i := \rho \circ \cdots \circ \rho$ ($i$ times). 
Then: 

(1)\ $\rho^i = \rho^{2r} $ if $i\ge 2r$; 

(2)\ $\hatE := \rho^{2r}(E)$ is a closed bundle; 

(3)\ $\dim \hatE_a $ is an upper-semicontinuous 
function of $a\in X$. 
\endproclaim

\demo{Proof}
Set 
$$ d_i (a) \ = \ \dim \rho^i (E)_a \ , \qquad 
\lambda_i (a) \ = \ \inf\limits_{\delta> 0} 
\sup\limits_{\sigma (a,x) < \delta} \ d_i (x) \ , $$ 
for all $a\in X$ and $i = 0$, $1,\dots, $ where
$\sigma (\cdot , \cdot)$ denotes the metric on $X$.
Then, for all $i$, $\lambda_i$ is upper-semicontinuous
and $d_i \le \lambda_i \le d_{i+1}$ (the latter inequality 
since $\overline{\rho^i (E)} \subset \rho^{i+1} (E)$).
Let 
$$ 
G_i \ := \ \rint \{ a \in X : \ d_i (a) = d_{i+1} (a) \} 
$$
and let $Z_i := X\bs G_i$.  Then $G_i \subset G_{i+1} $ for 
all $i$, since, for all $a\in G_i$, $\rho^i (E)_a = 
\rho^{i+1} (E)_a$ and therefore $\rho^{i+1} (E)_a = 
\rho^{i+2} (E)_a$ by locality (Definition 3.2, property 
(2)).  Thus $Z_i \supset Z_{i+1}$ for all $i$. 

We claim that $d_{i+2} (a) > i/2$, for all $a\in Z_i $, 
$i\in \IN$. First, this holds for $i=0$:  Otherwise, 
there exists $a\in Z_0$ such that $d_2 (a) = 0$. 
Then $\lambda_1 (a) = 0$, hence $d_0 (x) = d_1 (x) = 0$
in a neighbourhood of $a$, so that $a \in G_0$
(a contradiction).  The claim is true for $i=1$ because, 
if $a\in Z_1$, then $a\in Z_0 $ so that $d_3 (a) \ge 
d_2 (a) \ge 1$.  Consider $i \ge 2$ and suppose that 
$d_{j+2} (a) > j/2$ for all $a\in Z_j$, when $j<i$. 
Set $\Sigma_i := \{ a\in X : d_i (a) < d_{i+1} (a) \}$. 
Then $Z_i = \overline{\Sigma}_i$, so $\Sigma_i \subset Z_i \subset
Z_{i-2}$.  Therefore, if $a\in \Sigma_i$, then 
$d_i (a) > (i-2)/2$, so that $d_{i+1} (a) > i/2$.  
It follows that, for all $a\in Z_i$, $d_{i+2} (a) \ge 
\lambda_{i+1} (a) > i/2$.  This proves the claim, by induction.

It follows that $Z_{2r} = \emptyset$, so (1) holds. 
(2) follows because $\overline{\widehat{E}} \subset \rho (\hatE) \subset
\hatE$, and (3) because $d_{2r} (a) \le \lambda_{2r} (a) 
\le d_{2r+1} (a) $, $a\in X$. 
\endprf
\enddemo

\example{Example 3.4} \cite{G}. 
For any bundle $E \subset X \times V$, define 
$$
\tilE \ := \ \bigcup\limits_{a\in X} \{ a\} \times \Span E_a 
$$
where $\Span $ denotes the linear span).  
Set $\lambda (E) := \widetilde{\barE}$. 
Then  $\lambda$ is a Glaeser operation and 
$\lambda (E) \subset \rho (E)$ for any Glaeser operation 
$\rho$. 
\endexample

\definition{Definition 3.5} \cite{G}.
Let $X$ be a closed subset of $\IR^n$. 
Define 
$$
\align
{\text{ ptg}} (X) \ &:= \ \{ 
(a,\sigma u) \in X \times \IR^n : \  
\sigma \in \IR \ , \ u= \lim\limits_{j\to \infty}
\frac{x_j - y_j}{|x_j - y_j | } \ , \\
&\qquad\qquad\qquad {\text{ where }} (x_j), (y_j) \subset X , \  x_j \ne y_j , 
{\text{ for all }} j \ , \\
&\qquad\qquad\qquad {\text{ and}} \lim\limits_{j\to \infty} x_j = a = 
\lim\limits_{j\to \infty} y_j \} \ ; \\
\tau (X) \ &:= \ \widehat{\widetilde{{\text {ptg}}(X)}} \ ,\\
\endalign
$$
where the ``saturation'' $\widehat{\phantom{00}}$ is with 
respect to the Glaeser operation $\lambda$ of Example 3.4. 
We call $\tau (X)$ the {\it (linearized) paratangent
bundle} of $X$.  If $a\in X$, the fibre $\tau_a (X)$ 
is called the {\it paratangent space} of $X$ at $a$. 
\enddefinition

Let $f:X\to \IR$ be a continuous function. 
(We identify a function with its graph and therefore write
$\tau (f) = \tau (\text{graph}\, f)$.)  We can consider 
$\tau (f)$ as a bundle over $X$, so that $\tau (f) \subset
\tau (X) \times \IR$. 

\proclaim{Theorem 3.6} \cite{G}, \cite{Br}. 
Let $X$ be a closed subset of $\IR^n$ and let $f:X \to \IR$ 
be a continuous function.  Then there exists $F\in \cC^1 (\IR^n)$
such that $F|X = f$ if and only if $\tau (f) : \tau (X) \to \IR$
(i.e., $\tau (f)$ is the graph  of a function $\tau (X) \to \IR$). 
In this case, each $\tau_a (f) :\tau_a (X) \to \IR$, $a\in X$, 
is the restriction to the paratangent space $\tau_a (X)$ 
of the derivative of $F$. 
\endproclaim

\head 4. Higher-order paratangent spaces \endhead

\subhead The remainder term in Taylor's theorem \endsubhead 
Let $X\subset U \subset \IR^n$, where $U$ is open and $X$
is closed in $U$.  Let $p\in \IN$.  
Let $F= (F^\alpha)_{\alpha\in\IN^n , |\alpha|\le p }$, 
where each $F^\alpha : X \to \IR$.  If $a\in X$, define 
$$
(T^p_a  F) (x) \ := \ \sum\limits_{|\alpha|\le p} 
\frac{1}{\alpha !} F^\alpha (a) (x-a)^\alpha \ . 
$$ 
Let $D^\alpha := 
\partial^{|\alpha|}/\partial x_1^{\alpha_1} \cdots 
\partial x_n^{\alpha_n}$, $\alpha \in \IN^n$. 
If $|\alpha| \le p$ and $b\in X$, set 
$$
\align
(R^p_a F)^\alpha (b) 
\ &:= \ F^\alpha (b) - D^\alpha (T^p_a F) (b) \tag 4.1 \\
\ &= \ F^\alpha (b) - \sum\limits_{|\beta|\le p - |\alpha | }
\frac{1}{\beta!} F^{\alpha+ \beta} (a) 
(b-a)^\beta \ , \\
\delta_\alpha (a,b) 
\ &:= \ \frac{(R^p_a F)^\alpha (b)}{|b-a|^{p-|\alpha|}}\ . 
\tag 4.2 
\endalign
$$
We recall Whitney's extension theorem \cite{W1}: 

\proclaim{Theorem 4.3}
Let $F^\alpha : X\to \IR$, 
$|\alpha|\le p$.  Then there exists $f\in \cC^p (U)$ such 
that $(D^\alpha f )| X = F^\alpha$, for all 
$|\alpha | \le p$, if and only if $\delta_\alpha (a,b) \to 0$
if $a$, $b\in X$ and $|\alpha | \le p$, as $|a-b| \to 0$. 
\endproclaim 

We say that $F= (F^\alpha)$ is a $\cC^p$ 
{\it Whitney field} on $X$ if it satisfies the conditions 
of Theorem 4.3.  

Let $\cP_p = \cP_p (\IR^n)$ denote the real vector space of 
polynomial functions on $\IR^n$ of degree at most $p$.  Let 
$\xi \in \cP^*_p$, where $\cP^*_p$ denotes the dual of 
$\cP_p$. Set
$$ \xi (F,a) \ := \ \xi (T^p_a F) \ =
\ \sum\limits_{|\alpha|\le p} 
F^\alpha (a) \xi_\alpha (a) \ , 
\tag 4.4 $$ 
where 
$$ 
\xi_\alpha (a) \ := \ \xi \left( \frac{1}{\alpha!} 
(x-a)^\alpha \right) \ . 
\tag 4.5
$$ 

If $\eta \in \cP^*_p$ and $b\in X$, then 
$$\align
\eta (F,b) 
&= \sum\limits_{|\alpha | \le p} F^\alpha (b) \eta_\alpha  (b) \\
&= \sum\limits_{|\alpha | \le p} 
   \left( \sum\limits_{|\beta | \le p - |\alpha| } 
   \frac{1}{\beta!} F^{\alpha +\beta} (a) (b-a)^\beta + 
   \delta_\alpha (a,b) |b-a|^{p-|\alpha|} \right) 
   \eta_\alpha  (b) \\
&= \sum\limits_{|\alpha | \le p} \sum\limits_{\beta \le \alpha} 
   \frac{1}{\beta!} 
   F^\alpha (a) (b-a)^\beta \eta_{\alpha - \beta} (b) + 
   \sum\limits_{|\alpha | \le p} \delta_\alpha (a,b) | b-a|^{p-|\alpha|}
   \eta_\alpha (b) \\
&= \sum\limits_{|\alpha | \le p} \sum\limits_{\beta \le \alpha} 
   \frac{1}{\beta!} 
   F^\alpha (a) (b-a)^\beta \eta 
   \left( \frac{1}{(\alpha - \beta)!} (x-b)^{\alpha -\beta} \right)
   + \sum\limits_{|\alpha | \le p} \delta_\alpha (a,b) 
   |b-a|^{p-|\alpha|} \eta_\alpha (b) \\
&= \sum\limits_{|\alpha | \le p} \frac{1}{\alpha !} 
   F^\alpha (a) \eta 
   \left( \sum\limits_{\beta\le \alpha} {\alpha\choose \beta} 
   (b-a)^{\beta} (x-b)^{\alpha-\rho} \right)
   + \sum\limits_{|\alpha | \le p} \delta_\alpha (a,b) 
   |b-a|^{p-|\alpha|} \eta_\alpha (b) \\
&= \sum\limits_{|\alpha | \le p}  F^\alpha (a) \eta 
   \left( \frac{1}{\alpha!} (x-a)^\alpha \right) 
   + \sum\limits_{|\alpha | \le p} \delta_\alpha (a,b) 
   |b-a|^{p-|\alpha|} \eta_\alpha (b) \ . \\
\endalign
$$
(If $\alpha = (\alpha_1,\dots, \alpha_n)$ and 
$\beta = (\beta_1,\dots,\beta_n)$, then $\beta \le \alpha$
means $\beta_i \le \alpha_i $, $i= 1,\dots, n$.) 
Therefore, 
$$
\eta (F,b) \ = \ \eta (F,a) + \sum\limits_{|\alpha | \le p} 
\delta_\alpha (a,b) | b-a|^{p-|\alpha|} \eta_\alpha (b) ;
\tag 4.6 
$$ 
$$ 
\xi (F,a)+\eta (F,b) \ = \ (\xi + \eta) (F,a) + 
\sum\limits_{|\alpha | \le p} \delta_\alpha (a,b) 
|b-a|^{p-|\alpha|} \eta_\alpha (b) \ . 
\tag 4.7 
$$ 
We will use the following lemma only in the case $k=1$ 
(but see Remarks 4.13(2) and Final Remarks 5.5). 

\proclaim{Lemma 4.8}
Let $X \subset U \subset \IR^n$, where $U$ is open 
and $X$ is closed in $U$.  Let $(a_{ij}) = (a_{i1}, a_{i2},\dots)$
and $(\xi_{ij}) = (\xi_{i1} , \xi_{i2},\dots ) $ denote
sequences in $X$ and $\cP^*_p$, respectively, 
for $i = 0 , 1,\dots, k$, such that: 

(1)\ The sequences $(a_{ij})$, 
$i=0, 1,\dots, k$, converge to a common point $a\in X$, 
and $\sum_{i=0}^k \xi_{ij}$ converges to $\xi \in \cP^*_p$. 

(2)\ $|a_{ij} - a_{0j}|^{p-|\alpha|} 
|\xi_{ij,\alpha} (a_{ij}) | \le c $, for all 
$i$, $j$ and $|\alpha | \le p$ (where $c$ is a constant). 

If $F= (F^\alpha)_{|\alpha|\le p} $ is a $\cC^p$ Whitney field
on $X$, then 
$$ \xi (F,a) \ = \ \lim\limits_{j\to \infty} 
\sum\limits^k_{i=0} \xi_{ij} (F,a_{ij}) \ . $$
\endproclaim 

\demo{Proof}
For each $j=1,2,\dots,$
$$ \align 
& \xi (F,a) - \sum\limits^k_{i=0} \xi_{ij} (F,a_{ij}) \\
&= \xi (F,a) - \sum\limits^k_{i=0} \xi_{ij} (F,a ) + 
   \sum\limits^k_{i=0} \big(\xi_{ij} (F,a ) - \xi_{ij} (F,a_{0j})\big)
   + \sum\limits^k_{i=0} \big(\xi_{ij} (F, a_{0j}) - \xi_{ij}
   (F,a_{ij})\big)\\
&= \left( \xi - \sum\limits^k_{i=0} \xi_{ij}\right) (T^p_a F) + 
   \sum\limits^k_{i=0} \xi_{ij} (T^p_a F - T^p_{a_{0j}} F) 
   + \sum\limits^k_{i=0} \xi_{ij} (T^p_{a_{0j}} F - T^p_{a_{ij}} F)\\
&= \left( \xi - \sum\limits^k_{i=0} \xi_{ij}\right) (T^p_a F) +
   \sum\limits^k_{i=0} \sum\limits_{|\alpha|\le p} 
   \delta_\alpha (a_{0j}, a) |a-a_{0j}|^{p-|\alpha|} 
   \xi_{ij , \alpha } (a)\\
&{\hskip 2truein} -\sum\limits^k_{i=0} \sum\limits_{|\alpha|\le p} 
\delta_\alpha (a_{0j} , a_{ij}) | a_{ij} - a_{0j}|^{p-|\alpha|} 
\xi_{ij , \alpha } (a_{ij}) \ . \\
\endalign
$$ 
Each of the three terms tends to $0$ as $j\to \infty$. 
\endprf
\enddemo

\subhead The paratangent bundle of order $p$ \endsubhead
We consider Glaeser operations on bundles of subspaces of $\cP^*_p$. 
Let $X\subset U \subset \IR^n$, where $U$ is open and $X$ is closed 
in $U$. Let $E\subset X \times \cP^*_p$ be any bundle of linear 
subspaces of $\cP^*_p$ over $X$. Define 
$$ 
\align
\Delta E \ := \ \{ (a , b , \xi+\eta) : \ &a,b\in X ,\ \xi \in E_a ,\
\eta \in E_b , \tag 4.9\\
&|a-b|^{p-|\alpha|} |\eta_\alpha (b)| \le 1 \ 
{\text{ for all }} |\alpha | \le p \}\ . 
\endalign
$$
Let $\pi : X^2 \times \cP^*_p \to X\times \cP^*_p$ denote the
projection $\pi (a,b,\xi) = (a,\xi) $. Define 
$$ 
E' \ := \ \pi (\overline{\Delta E} \cap \{ (a,a,\xi) :
\ a\in X , \xi \in \cP^*_p \} ) \ . 
\tag 4.10 
$$
Clearly, $\barE \subset E'$.  We define a Glaeser 
operation 
$$ 
\rho (E) \ := \ \widetilde{E'} \ := \ \bigcup\limits_{a\in X} 
\{ a \} \times \Span E'_a \ . 
\tag 4.11
$$ 
(Recall Example 3.4.)

\definition{Definition 4.12}
Let $X\subset U \subset \IR^n$ be as above, and set 
$$ 
E \ := \ \{ (a,\lambda \delta_a ) : \ a\in X , \lambda \in \IR \}\ , 
$$
where $\delta_a \in \cP^*_p$ denotes the delta-function 
$\delta_a (P) := P (a)$, $P\in \cP_p$.  Define 
$$ \tau^p (X) \ := \ \hatE \ , $$ 
where $\hatE$ denotes the saturation of $E$ with 
respect to the Glaeser operation (4.11) (cf. Lemma 3.3). 
We call $\tau^p (X)$ the {\it (linearized) paratangent 
bundle} of $X$ of {\it order} $p$. 
\enddefinition

\remark{Remarks 4.13}
(1) Consider $q \geq p$. Recall that $X \times \cP_p^*$ embeds
in $X \times \cP_q^*$ as a closed subbundle, where, for each 
$a \in X$, the fibre $\cP_p^*$ over $a$ is identified as in 
Section 2 with $(\undm^{p+1}_a)^\perp \subset \cP_q^*$. If
$b \in X$ and $\eta \in (\undm^{p+1}_b)^\perp$, then 
$\eta_\alpha (b) = 0$ when $p < |\alpha | \leq q$.
It follows from the definition above that $\tau^p (X)$ is a
closed subbundle of $\tau^q (X)$.
\smallskip

(2) The definition above involves distributions with values
in $\cP^*_p$ supported at pairs of points $a, b \in X$, according
to (4.9), and suffices for all results in this paper. But a more
general definition of $\tau^p (X)$ involving distributions
supported at $k+1$ points (where $k \geq p$) is necessary for
our main conjecture in Section 1. See Final Remarks 5.5. We
have stated Lemma 4.8 and Lemma 5.1 below for distributions
supported at $k+1$ points in order that they be available more
generally.
\endremark
\smallskip

Now let $\Phi \subset X \times (\cP^*_p \times \IR)$ be a bundle 
of linear subspaces of $\cP^*_p \times \IR$ over $X$.  Define 
$$
\align 
\Delta \Phi \ &:= \ \{ (a,b, \xi +\eta , \lambda + \mu ) : \ a,b\in X ,\ 
(\xi ,\lambda )\in \Phi_a\ ,\ (\eta,\mu) \in \Phi_b \ , \tag 4.14 \\
&\qquad\qquad\qquad\qquad\qquad |a-b|^{p-|\alpha|} | \eta_\alpha (b) | \le 1 
{\text{ for all }} |\alpha| \le p \} \ ; \\
\Phi' \ &:= \ \pi (\overline{\Delta \Phi} \cap \{ (a,a,\xi,\lambda): 
\ a\in X , \xi \in \cP^*_p \ , \lambda \in \IR\} ) \ , 
\tag 4.15 
\endalign
$$ 
where $\pi :X^2 \times \cP^*_p \times \IR \to X \times 
\cP^*_p \times \IR $ is the projection $\pi (a,b,\xi , \lambda) = 
(a,\xi ,\lambda )$.  As before, $\overline{\Phi} \subset \Phi'$. 

\definition{Definition 4.16}
Define 
$$ 
\nabla^p f \ := \ \widehat\Phi \ , 
$$
where 
$$
\Phi \ := \ \{ (a,\lambda \delta_a , \lambda f(a)) : \ a\in X , 
\lambda \in \IR \} $$
and $\widehat\Phi $ denotes the saturation with respect to 
the Glaeser operation $\rho (\Phi ) := \widetilde{\Phi'}$. 
\enddefinition

Clearly, $\nabla^p f \subset \tau^p (X) \times \IR$. 
Theorem 1.1 is a restatement of Theorem 4.18 below.

\proclaim{Lemma 4.17}
Let $X\subset U \subset \IR^n$, where $U$ is open and $X$ 
is closed in $U$.  Let $f : X \to \IR$.  Let $p\in \IN$
and suppose there is a $\cC^p$ Whitney field 
$F = (F^\alpha )_{|\alpha| \le p}$ on $X$ such that $F^0 = f$. 
Consider the bundles $E$, $\Phi$ over $X$ and the Glaeser 
operations $\rho$ of Definitions 4.12 and 4.16. Then, 
for each $i\in \IN$, 
$$ 
\rho^i (\Phi ) : \ \rho^i (E) \to \IR \ ; 
$$
moreover, if $a\in X$ and $\xi \in \rho^i (E)_a $, then 
$$ 
\rho^i (\Phi) (\xi ) \ = \ \xi (T^p_a F) \ = \ \xi (F,a) \ . 
$$ 
\endproclaim 

\demo{Proof}
First consider $i=1$. Any element of $\Phi'$ can be expressed
$$ 
\lim\limits_{j\to \infty} 
(a_{0j}, a_{1j}, \xi_{0j} + \xi_{1j} , 
\lambda_{0j} f (a_{0j}) + \lambda_{1j} f (a_{1j})) \ , $$
where $\xi_{ij} = \lambda_{ij} \delta_{a_{ij}} $, 
$i=0,1$, $j=1,2,\dots,$ and $(a_{ij})$, $(\xi_{ij})$
satisfy the hypotheses of Lemma 4.8 (case $k=1$). 
Of course , $\lambda_{ij} f (a_{ij}) = \xi_{ij} (F,a_{ij})$, 
for all $i$, $j$.  By Lemma 4.8, 
$$ 
\lim\limits_{j\to \infty} 
(\lambda_{0j} f (a_{0j}) + \lambda_{1j} f (a_{1j})) \ = 
\ \lim\limits_{j\to \infty} 
(\xi_{0j} + \xi_{1j}) (F,a) \ . 
$$ 
Therefore, $\Phi' : E' \to \IR$ and, for all
$\xi \in E'_a $, $a\in X$, 
$\Phi' (\xi) = \xi (F,a) = \xi (T^p_a F)$. 
It follows that $\rho (\Phi) : \rho (E) \to \IR$ and, 
for all $\xi \in \rho (E)_a $, 
$\rho (\Phi) (\xi) = \xi (T^p_a F) = \xi (F,a) $. 

The result then follows from Lemma 4.8, by induction on $i$.
\endprf
\enddemo

\proclaim{Theorem 4.18}
Let $X\subset U \subset\IR^n$, where $U$ is open and $X$ is 
closed in $U$.  Let $f:X\to \IR$. Let $p\in \IN$ and suppose 
there is a $\cC^p$ Whitney field $F= (F^\alpha)_{|\alpha| \le p}$
on $X$ such that $F^0 = f$.  Then 
$$ \nabla^p f : \ \tau^p (X) \to \IR \ ; $$ 
moreover, if $a\in X$ and $\xi \in \tau^p_a (X) \subset \cP^*_p$, 
then 
$$ 
\nabla^p f(\xi) \ = \ \xi (T^p_a F) \ = \ \xi (F,a) \ . 
$$ 
\endproclaim

This is an immediate consequence of Lemma 4.17.

\remark{Remarks 4.19}
Let $f:X \to \IR$.  If $\nabla^p f: \tau^p (X) \to \IR$, 
then $\nabla^p f $ (as well as $f$) is necessarily 
continuous and, for all $a\in X$, the induced function on the 
fibre $\nabla^p_a f : \tau^p_a (X) \to \IR$ is necessarily 
linear. Consider $q \geq p$. Then $\nabla^p f $ is a closed
subbundle of $\nabla^q f $ (cf. Remarks 4.13(1)); if 
$\nabla^q f : \tau^q (X) \to \IR$, then 
$\nabla^p f : \tau^p (X) \to \IR$ and $\nabla^p f$ is the
restriction of $\nabla^q f$.
\endremark
\smallskip

We will show that Theorem 1.2 is a consequence of Theorem 4.21 
below. 

\proclaim{Lemma 4.20}
Suppose that $X\subset U \subset \IR^n$, where $U$ is open and 
$X$ is the closure of $\, \text{\rm int} \, X$ in $U$.  Then 
$\tau^p (X) = X\times \cP^*_p$, for all $p\in \IN$. 
\endproclaim 

\demo{Proof}
It is enough to show that $\tau^p_a (X) = \cP^*_p$, 
where $a\in \rint X$. If $\alpha \in \IN^n $, $|\alpha|\le p$, 
and $b\in X$, define $D^\alpha (b) \in \cP^*_p $ by 
$$ D^\alpha (b) (P) \ := \ (D^\alpha P) (b) \ , \quad 
P\in \cP_p \ . $$
If $|\alpha| < p$, then 
$$ 
\frac{D^\alpha (b) - D^\alpha (a)}{|b-a|} 
\ \longrightarrow\ 
\sum\limits_{i=1}^n u_i D^{\alpha + (i)} (a) 
$$
if $b\to a$ and $\displaystyle{\frac{b-a}{|b-a|} 
\to u = (u_1,\dots,u_n)}$, where 
$(i)$ denotes the multiindex with $1$ in the $i$'th place 
and $0$ elsewhere.  Let 
$\displaystyle {\eta = \frac{D^\alpha (b)}{|b-a|}}$
(where $|\alpha| < p$).  Then, for all $\gamma \in \IN^n$, 
$|\gamma | \le p$, 
$$\align
|a-b|^{p-|\gamma|} \eta_\gamma (b) \ &=
\ |a-b|^{p-|\gamma|-1} D^\alpha 
\left( \frac{1}{\gamma!} (x-b)^\gamma \right) (b)\\
\ &= \ \cases 
                  0  &,\quad \gamma \ne \alpha \\
|a-b|^{p-|\alpha|-1} &,\quad \gamma = \alpha \ . \endcases \\
\endalign
$$
By Definition 4.12, it follows by induction on $|\gamma|$ 
that $D^\gamma (a) \in \tau^p_a (X)$, for all 
$|\gamma | \le p $; i.e., 
$\tau^p_a (X) = \cP^*_p $. 
\endprf
\enddemo

\proclaim{Theorem 4.21}
Let $X\subset U\subset \IR^n$, where $U$ is open 
and $X$ is the closure of $\, \text{\rm int} X$ in $U$. Let $f:X\to \IR$.
Let $p\in \IN$. Suppose that 
$$\nabla^p f : \ \tau^p (X) \to \IR \ . $$
Then there is a $\cC^p$ Whitney field $F = (F^\alpha)_{|\alpha|\le p}$
on $X$ such that $F^0 = f$. 
\endproclaim 

\demo{Proof}
By Lemma 4.20, we have 
$$ \nabla^p f : \ X \times \cP^*_p \to \IR \ . $$ 
Define 
$$ F^\alpha (a) \ := \ (\nabla^p f) (a,D^\alpha (a)) \ , \quad
a\in X \ , $$
for all $\alpha \in \IN^n$, $|\alpha | \le p$. 
(We use the notation of the proof of Lemma 4.20.)  Let $c\in X$
and $\alpha \in \IN^n$, $|\alpha| \le p$. Then 
$$ 
\frac{D^\alpha (b) - \sum\limits_{|\beta|\le p - |\alpha|}
\dfrac{1}{\beta!} (b-a)^\beta D^{\alpha+\rho} (a)}
{|b-a|^{p-|\alpha|}} 
\ \longrightarrow \ 0 
$$ 
as $a$, $b\to c$, where $a$, $b\in X$, $a\ne b$. 
(In fact, this element of $\cP^*_p$ equals zero since, 
for all $P \in \cP_p$, $(D^\alpha P) (b) = 
(T^{p-|\alpha|}_a D^\alpha P ) (b)$.) If 
$\displaystyle{\eta := \frac{D^\alpha(b)}{|b-a|^{p-|\alpha|}}}$, then 
$$ 
|a-b|^{p-|\gamma|} \eta_\gamma (b) \ = 
\ \cases
0 &,\quad \gamma\ne \alpha \\
1 &,\quad \gamma = \alpha \ . \endcases
$$
Hence
$$ 
\frac{(\nabla^p f)(b,D^\alpha (b)) -\sum\limits_{|\beta|\le p-|\alpha|}
\dfrac{1}{\beta!} (b-a)^\beta (\nabla^p f) (a , D^{\alpha+\beta} (a))}
{|b-a|^{p-|\alpha|}} 
\ \longrightarrow \ 0 
$$ 
as $a$, $b\to c$ in $X$, $a\ne b$; in other words, 
$$ 
F^\alpha (b) - \sum\limits_{|\beta| \le p -|\alpha|} 
\frac{1}{\beta!} F^{\alpha +\beta} (a) 
(b-a)^\beta \ = \ o (|b-a|^{p-|\alpha|}) \ , 
$$ 
as required. 
\endprf
\enddemo

\remark{Remark 4.22}
The following is a simple generalization of Theorem 4.21: 
Consider $1\le m\le n$ and 
$X\subset U \times \{ 0 \} \subset \IR^m \times \IR^{n-m} 
= \IR^n$, where $U$ is open in $\IR^m$ and $X = \overline{\rint X}$
as a subset of $U$.  Let $f: X\to \IR$.  Suppose that 
$\nabla^p f : \tau^p (X) \to \IR$.  Then there is a $\cC^p$ 
Whitney field $F = (F^\alpha)_{\alpha \in \IN^n,|\alpha | \le p}$
on $X$ such that $F^0 = f $ and $F^\alpha=0$ when 
$\alpha_{m+1} + \cdots + \alpha_n > 0 $. 
\endremark
\smallskip

\definition{Definition 4.23}
Let $X$, $Y$ denote metric spaces and $V$, $W$ finite-dimensional
real vector spaces.  Let $E\subset X \times V$ and 
$F \subset Y \times W$ be bundles (of linear subspaces 
of $V$ and $W$, respectively).  A {\it morphism} $E\to F$
is a continuous mapping $\psi: E \to F$ of the form 
$\psi (a,v) = (\varphi(a) , \psi_1 (a,v))$, where $(a,v) \in E$, 
such that, for all $a\in X$, $\psi_1 (a,\cdot) : E_a \to 
F_{\varphi(a)}$ is linear.  An {\it isomorphism} is a 
morphism with a continuous inverse (which is necessarily a 
morphism). 
\enddefinition

Suppose that $U_1$, $U_2$ are open subsets of $\IR^n$ and that 
$X_1$, $X_2$ are closed subsets of $U_1$, $U_2$, respectively.
Let $\sigma: U_1 \to U_2$ be a $\cC^p$ diffeomorphism
$(p \in \IN)$ such that $\sigma (X_1) = X_2$.  Clearly, 
$\sigma $ induces an isomorphism $\sigma_* : \tau^p (X_1)\to 
\tau^p (X_2)$. (See Theorem 5.2 below.)
If $f_1 : X_1 \to \IR$ and $f_2:X_2 \to \IR$
are functions such that $f_1 = f_2 \circ \sigma$, then 
$\sigma $ induces an isomorphism $\sigma^* : \nabla^p f_1 \to 
\nabla^p f_2 $. These observations can be used to generalize
the results above to manifolds.  Theorem 1.2 is a special 
case of the following.

\proclaim{Theorem 4.24}
Let $X\subset M \subset U \subset \IR$, where $U$ is open, 
$M$ is a closed $\cC^p$ submanifold of $U$, and 
$X= \overline{\text{\rm int} X}$ as a subset of $M$.  Let $f:X\to \IR$.
If $\nabla^p f: \tau^p (X) \to \IR$, then $f$ is the 
restriction of an element of $\cC^p (U)$. 
\endproclaim

\head 5. Composite functions\endhead

Let $U$, $V$ be open subsets of $\IR^n$, $\IR^m$ (respectively)
and let $\varphi : V \to U$ be a $\cC^p$ mapping.  Let 
$b\in V$, $a = \varphi (b)$.  Then $\varphi$ induces a linear 
mapping 
$$ 
\align
\varphi^*_b : \ \cP_p (\IR^n) &\to \cP_p (\IR^m) \\
P &\mapsto T^p_b (P \circ \varphi) \ ; \\
\endalign
$$
i.e., if $x = (x_1,\dots, x_n)$, $y= (y_1,\dots, y_m)$ denote 
the coordinates of $\IR^n$, $\IR^m$ (respectively) and 
$\varphi = (\varphi_1,\dots, \varphi_n)$, then 
$\varphi^*_b (P)$ is given by substituting $T^p_b \varphi = 
(T^p_b \varphi_1,\dots, T^p_b \varphi_n)$ into 
$P(x) = (T^p_a P) (x)$ and truncating terms involving 
$(y-b)^\beta$ where $|\beta | > p$.  By duality, there is a 
linear mapping 
$$ \varphi_{*b} : \ \cP_p (\IR^m)^* \to \cP_p (\IR^n)^* \ ; $$
i.e., $\varphi_{*b} (\eta) (P) = \eta (\varphi^*_b (P))$, 
where $\eta \in \cP_p (\IR^m)^*$ and $P\in \cP_p (\IR^n)$. 

Note that $\varphi_{*b} (\delta_b) = \delta_a$. 
We will need the following lemma only in the case $k=1$
(cf. Lemma 4.8 and Remarks 4.13(2)).

\proclaim{Lemma 5.1}
Let $X$, $Y$ be closed subsets of $U$, $V$ (respectively), 
where $U\subset \IR^n$, $V\subset \IR^m$ are open. 
Let $\varphi : V\to U$ be a $\cC^p$ mapping such that 
$\varphi (Y) \subset X$.  Let $(b_{ij}) = (b_{i1} ,b_{i2},\dots )$
and $(\eta_{ij}) = (\eta_{i1}, \eta_{i2},\dots )$ denote
sequences in $Y$ and $\cP_p (\IR^m)^*$, respectively, for 
$i=0,1,\dots, k$ such that: 

(1)\ The sequences $(b_{ij})$, $i=0,1,\dots,k$, 
converge to a common point $b\in Y$, and 
$\left( \sum\limits_{i=0}^k \eta_{ij} \right)$ converges to 
$\eta\in \cP_p (\IR^m)^*$. 

(2)\ $|b_{ij} - b_{0j}|^{p-|\beta|} | \eta_{ij,\beta}
(b_{ij}) | \le c$, for all $i$, $j$ and $\beta \in \IN^m$, 
$|\beta | \le p$, where $c$ is a constant. 

Set $a_{ij} = \varphi (b_{ij}) \in X$ and 
$\xi_{ij} = \varphi_{*b_{ij}} (\eta_{ij})$, 
for all $i$, $j$, and set $a=\varphi (b)$, 
$\xi = \varphi_{*b} (\eta)$. 
Then: 

(1$'$)\ $(a_{ij})$ converges to $a\in X$, for all 
$i=0,1,\dots , k$, and 
$\left( \sum\limits_{i=0}^k \xi_{ij} \right)$ converges to $\xi$. 

(2$'$)\ $| a_{ij} - a_{0j}|^{p-|\alpha|} |
\xi_{ij ,\alpha} (a_{ij})| \le c'$, for all $i$, $j$ and 
$\alpha  \in \IN^n$, $|\alpha | \le p$, where $c'$ is a 
constant. 
\endproclaim 

\demo{Proof}
Obviously, each $(a_{ij})$ converges to $a$.  Let $P \in \cP_p (\IR^n)$.
Then, for each $j$, 
$$ \align
\xi (P) - \sum\limits_{i=0}^k \xi_{ij} (P) 
\ &= \ \eta (\varphi^*_b (P)) - \sum\limits_{i=0}^k \eta_{ij}
(\varphi^*_{b_{ij}} (P)) \\
\ &= \ \eta (G,b) - \sum\limits_{i=0}^k \eta_{ij} (G,b_{ij}) \ , \\
\endalign
$$
where $G$ denotes the $\cC^p$ Whitney field on $Y$ induced by 
$P \circ \varphi$.  Therefore $(1')$ follows from Lemma 4.8. 

There is a constant $C$ such that, for all $i$, $j$ and 
$\alpha \in \IN^n$, $|\alpha | \le p$, 
$$ 
\align
|a_{ij} &- a_{0j} |^{p-|\alpha|} | \xi_{ij,\alpha} (a_{ij})| \\
&\le \ C|b_{ij}-b_{0j}|^{p-|\alpha|} 
\left| \eta_{ij} \left(\varphi^*_{b_{ij}} \left( \frac{1}{\alpha!} 
(x-a_{ij})^\alpha\right) \right) \right| \\
&\le \ C |b_{ij}-b_{0j}|^{p-|\alpha|} 
\sum\limits_{\beta\in \IN^m\atop |\alpha|\le |\beta|\le p}
\Lambda_\beta 
\Big( \big(D^\gamma \varphi (b_{ij})\big)_{1\le |\gamma| \le p}\Big)
\big| \eta_{ij,\beta} (b_{ij}) \big| \ , \\
\endalign
$$
where each $\Lambda_\beta$ is a polynomial function. 
Therefore, (2$'$) follows from (2). 
\endprf
\enddemo

\proclaim{Theorem 5.2}
Suppose $X\subset U \subset \IR^n$ and $Y\subset V \subset \IR^m$, 
where $U$, $V$ are open and $X$, $Y$ are closed in $U$, $V$
(respectively).  Let $\varphi : V \to U$ be a $\cC^p$ mapping 
such that $\varphi (Y) \subset X$.  Then: 

(1)\ $\varphi $ induces a bundle morphism 
$$ 
\varphi_* : \ \tau^p (Y) \to \tau^p (X)
$$ 
such that , if $b\in Y$ and $\eta \in \tau^p_b (Y) \subset 
\cP_p (\IR^m)^*$, then 
$$ 
\varphi_* (\eta) \ = \ \varphi_{*b} (\eta) \ . 
$$ 

Moreover, let $f:X\to \IR$ and set $g = f\circ \varphi: Y\to \IR$. 
Suppose that $\nabla^p f : \tau^p (X) \to \IR$.  Then: 

(2)\ $\nabla^p g : \tau^p (Y) \to \IR$ and, if $b\in Y$ and 
$\eta \in \tau^p_b (Y)$, then 
$$ \nabla^p g(\eta) \ = \ \nabla^p f (\varphi_* \eta)\ . $$ 

(3)\ Let $b\in Y$ and $a = \varphi (b)$.  Choose $P\in \cP_p 
(\IR^n)$ such that $P | \tau^p_a (X) : \tau^p_a (X) \to \IR$
coincides with $\nabla^p_a f$ (where we have identified
$\cP_p (\IR^n)$ with $\cP_p (\IR^n)^{**}$). 
Then, for all $\eta \in \tau^p_b (Y) \subset \cP_p (\IR^m)^*$, 
$$
\nabla^p g (\eta ) \ = \ \eta \big(\varphi^*_b (P)\big) \ . $$ 
\endproclaim 

\demo{Proof}
(1) follows from Lemma 5.1 and the definition of the paratangent
bundle in the same way that Lemma 4.17 is proved above using 
Lemma 4.8.

Consider $b_{ij} \in Y $ and $\eta_{ij} \in \tau^p_{b_{ij}} (Y)
\subset \cP_p (\IR^m)^*$, $i=0,1$, $j=1,2,\dots$, satisfying 
the hypotheses of Lemma 5.1 (case $k=1$).  Let $a_{ij}$, 
$\xi_{ij}$, $a$ and $\xi$ be as in Lemma 5.1.  Then, 
by Lemma 5.1 and Remark 4.19, 
$$ \align
\nabla^p f (\varphi_{*b} (\eta)) 
\ &= \ \nabla^p f(\xi) \\
&= \ \lim\limits_{j\to\infty}\big( \nabla^p f (\xi_{0j}) + 
\nabla^p f (\xi_{1j}) \big)\\
&= \ \lim\limits_{j\to\infty}\big( \nabla^p f 
\big(\varphi_{*b_{0j}} (\eta_{0j})\big) + \nabla^p f 
\big(\varphi_{*b_{1j}} (\eta_{1j})\big) \big)\ .  \\
\endalign
$$ 
(2) then follows in the same way that (1) is proved. 

To prove (3):  Let $\eta \in \tau^p_b (Y)$ and let 
$\xi = \varphi_{*b} (\eta) \in \tau^p_a (X) \subset \cP_p (\IR^n)^*$.
Then $\nabla^p f(\xi) = \xi (P)$, by the choice of $P$, so that 
$$ \align
\nabla^p f (\varphi_{*b} (\eta)) 
\ &= \ \nabla^p f(\xi) \\
&= \ \varphi_{*b} (\eta) (P) \\
&= \ \eta \big(\varphi^*_b (P)\big)\ , \\
\endalign
$$ 
and the result follows from (2).
\endprf
\enddemo

\proclaim{Corollary 5.3}
Suppose $X\subset U \subset \IR^n$ and $V\subset \IR^m$, 
where $U$, $V$ are open and $X$ is closed in $U$. 
Let $\varphi : V\to U$ be a $\cC^p$ mapping such that 
$\varphi(V) \subset X$.  Let $f: X \to \IR$ and set 
$g = f\circ\varphi$.  Suppose  that $\nabla^p f : \tau^p (X) \to 
\IR$.  Then: 

(1)\ $g\in \cC^p (V)$.

(2)\ $g$ is formally a composite with $\varphi$; i.e., for all 
$a\in X$, there exists $P\in \cP_p (\IR^n)$ such that 
$g-P \circ \varphi$ is $p$-flat at every point 
$b \in \varphi^{-1} (a)$. 
\endproclaim 

\demo{Proof}
(1) follows from Theorem 5.2 (2) and Theorem 4.21. 

Let $a\in X$. Choose $P \in \cP_p (\IR^n)$ as in Theorem 5.2 (3). 
Let $b\in \varphi^{-1} (a)$.  We will show that 
$$ 
T^p_b g \ = \ T^p_b (P \circ \varphi) \ . 
\tag 5.4
$$ 
Since $T^p_b (P\circ \varphi) = \varphi^*_b (P)$, (5.4) means 
that, for all $\eta \in \tau^p_b (V) = \cP_p (\IR^m)^*$, 
$$ 
\eta (T^p_b g) \ = \ \eta \big(\varphi^*_b (P)\big)\ . 
$$
But $\eta (T^p_b g) = (\nabla^p g) (\eta)$, by Theorem 4.18, 
so the result follows from Theorem 5.2 (3). 
\endprf
\enddemo

\subhead Differentiable functions on closed subanalytic sets
\endsubhead

\demo{Proof of Theorem 1.3}
Let $X$ be a compact subanalytic subset of $\IR^n$. By 
\cite{BM1, Thm. 0.1}, there is a compact real analytic
manifold $M$ and a real analytic mapping $\varphi : M \to \IR^n$ 
such that $\varphi(M) = X$.  By \cite{BMP}, there is a function 
$q= q_\varphi (p)$ from $\IN$ to itself such that, if $g \in \cC^q (M)$
and $g$ is formally a composite with $\varphi$, then there exists
$F \in \cC^p (\IR^n)$ such that $g = F \circ \varphi$; moreover,
if $S$ is a finite subset of $X$ and $g$ is $q$-flat on 
$\varphi^{-1}(S)$, then there exists $F$ with the additional
property that $F$ is $p$-flat on $S$.

Let $f:X\to \IR$.  Let $p\in \IN$ and suppose that 
$\nabla^q f : \tau^q (X) \to \IR$, where $q=q_\varphi (p)$. 
Let $g = f\circ \varphi$.  By Corollary 5.3 
(generalized to a manifold $V$), $g\in \cC^q (M)$ and $g$ 
is formally a composite with $\varphi$.  Therefore, 
$f\in \cC^p (X)$. 
\endprf
\enddemo

\remark{Remark 5.4}
If $X$ is a closed subanalytic subset of $\IR^n$, then 
$\tau^p (X)$ is a closed subanalytic subset of $\IR^n \times 
\cP_p (\IR^n)^*$. 
\endremark

\remark{Final Remarks 5.5}
(1) Let $X\subset U \subset \IR^n$, where $U$ is open and $X$ 
is closed in $U$.  Let $f:X\to \IR$.  Our definitions of 
$\tau^p (X)$ and $\nabla^p f $ involve limits of distributions
with values in $\cP_p (\IR^n)^*$ supported at two points.
We can generalize the definitions (and all constructions in the 
article) by using distributions supported at $k+1$ points, for any 
$k=1,2,\dots$.  We simply modify (4.9) and (4.10) in the following 
way: Let $E\subset X \times \cP_p (\IR^n)^*$ be any bundle 
of linear subspaces of $\cP_p (\IR^n)^* $ over $X$. 
Define 
$$
\align
\Delta_{k+1} E 
\ &:= \ \{ (a_0,a_1,\dots,a_k, \xi_0 + \xi_1 + \cdots + \xi_k ): \ 
a_i \in X \ , \ 
\xi_i \in E_{a_i} \ , \tag 5.6 \\
&\qquad \ |a_i - a_0|^{p-|\alpha|} | 
\xi_{i\alpha} (a_i) | \le 1 \ , \  
{\text{for all }} |\alpha | \le p \ , \ i=0,\dots,k \} \ ; \\
E'_{k+1} \ &:= \ \pi ( \overline{\Delta_{k+1} E} \cap 
\{ (a,a,\dots,a,\xi) : \ a\in X \ , \ \xi\in \cP_p 
(\IR^n)^* ) \ , 
\tag 5.7
\endalign
$$
where $\pi $ is the projection $\pi (a_0 , a_1 ,\dots, a_k ,\xi) 
= (a_0 , \xi)$.  We define $\tau^p_{k+1} (X) = \hatE$
as before, and $\nabla^p_{k+1} f \subset 
\tau^p_{k+1} (X) \times \IR$ also in a similar way. 

Of course $\tau^p_k (X) \subset \tau^p_{k+1} (X)$ and 
$\nabla^p_k f \subset \nabla^p_{k+1} f$, for all 
$k\ge 2$; in particular, if $\nabla^p_{k+1} f: \tau^p_{k+1} (X) 
\to \IR$, then $\nabla^p_k f: \tau^p_k (X) \to \IR$.
We have used only $\tau^p (X) = \tau^p_2 (X)$
in this article because it suffices for all the results. 
Our main conjecture in Section 1 should be
understood as requiring $\tau^p_{k+1} (X)$, where $k \geq p$. 
(For example, if
$$
X \ = \ \bigcup_{i=0}^p \{ (x,y) \in \IR^2: \ y = ix^2 \} \ ,
$$
then $k = p$ is necessary and sufficient.)

\demo{Questions}
Does there exist $r= r(X,p) \in \IN$ such that $\tau^p_k (X) 
= \tau^p_r (X)$ if $k\ge r$?  If $X$ is subanalytic, can we take
$r= p+1$?
\enddemo

(2) It is not difficult to see that the definition of 
$\Delta_{k+1} E$ above is equivalent to that given by 
replacing the condition
$$ 
|a_i - a_0 |^{p-|\alpha|} | \xi_{i\alpha} (a_i) |
\ \le \ 1 \ ,\quad i=0,\dots,k \ , 
$$ 
by the condition
$$ 
|a_i - a_0 |^{p-|\alpha|} | \xi_{i\alpha} (a_0) |
\ \le \ 1 \ ,\quad i=0,\dots,k \ . 
$$ 
(Likewise in Lemma 4.8).  It is not possible, however, to define 
$\tau^p_{k+1} (X)$ using limits
$$ \xi \ = \ \lim\limits_{j\to \infty} 
\sum\limits_{i=0}^k \xi_{ij} $$
(in the notation of Lemma 4.8) where condition (2) of Lemma 4.8
is replaced by the symmetric condition 
$$ 
|a_{ij} - a|^{p-|\alpha|} | \xi_{ij,\alpha } 
(a_{ij}) | \ \le \ c \ , 
\tag 5.8 
$$
for all $i$, $j$ and $|\alpha | \le p$. 

For example, let $(x_1,y_1) = (1,1)$ and, for each 
$j=1,2,\dots,$ define $(x_{j+1} ,y_{j+1})$ inductively as 
follows:  If $j$ is odd (respectively, even), let 
$(x_{j+1},y_{j+1})$ be the intersection point of the line 
through $(x_j , y_j)$ with slope $2$ (respectively, $-2$) and the 
arc $y=-x^2$, $x>0$ (respectively $y=x^2$, $x>0$).
Let $X = \{ 0 \} \cup \{ x_j: j \ge 1 \} \subset \IR$.
Define $F^0 (0) = 0$, $F^0 (x_j) = y_j$, for all $j$, and 
$F^1 (a) = 0$, for all $a\in X$.  Then 
$$ 
\lim\limits_{j\to \infty}
\frac{(R^1_0 F)^0 (x_j)}{|x_j - 0|} \ = \ 0 \ , 
$$
but 
$$
\frac{(R^1_0 F)^0 (x_j)}{|x_j - x_{j+1}|} 
$$
does not tend to zero as $j\to \infty$, so that $F$ 
is not a Whitney field.  Take $a_{0j} = x_j$, 
$a_{1j} = x_{j+1} $, 
$\displaystyle{\xi_{0j} = \frac{\delta_{x_j}}{x_j - x_{j+1}}}$ and 
$\displaystyle{\xi_{1j} = \frac{\delta_{x_{j+1}}}{x_j - x_{j+1}}}$,
for all $j$.  Then the condition (2) of Lemma 4.8 
(case $k=1$) is satisfied, but not the symmetric 
condition (5.8). 
\endremark

\enddocument